\documentclass[11pt]{article}
\usepackage{amsmath}
\usepackage{amssymb}
\usepackage{amsthm}
\usepackage{mathrsfs}
\begin{document}
\title{The intrinsic square function characterizations of weighted Hardy spaces}
\author{Hua Wang,\quad Heping Liu\,\footnote{E-mail address: wanghua@pku.edu.cn, hpliu@pku.edu.cn.}\\[3mm]
\footnotesize{School of Mathematical Sciences, Peking University, Beijing 100871, China}}
\date{}
\maketitle
\begin{abstract}
In this paper, we will study the boundedness of intrinsic square
functions on the weighted Hardy spaces $H^p(w)$ for $0<p<1$, where
$w$ is a Muckenhoupt's weight function. We will also give some
intrinsic square function characterizations of weighted Hardy spaces
$H^p(w)$ for $0<p<1$.\\[2mm]
\end{abstract}
\textbf{\large{1. Introduction and preliminaries}}
\par
First, let's recall some standard definitions and notations. The classical $A_p$ weight theory was first introduced by Muckenhoupt in the study of weighted
$L^p$ boundedness of Hardy-Littlewood maximal functions in [8].
Let $w$ be a nonnegative, locally integrable function defined on $\mathbb R^n$, all cubes are assumed to have their sides parallel to the coordinate axes.
We say that $w\in A_p$, $1<p<\infty$, if
$$\left(\frac1{|Q|}\int_Q w(x)\,dx\right)\left(\frac1{|Q|}\int_Q w(x)^{-\frac{1}{p-1}}\,dx\right)^{p-1}\le C \quad\mbox{for every cube}\; Q\subseteq \mathbb
R^n,$$
where $C$ is a positive constant which is independent of the choice of $Q$.\\
For the case $p=1$, $w\in A_1$, if
$$\frac1{|Q|}\int_Q w(x)\,dx\le C\,\underset{x\in Q}{\mbox{ess\,inf}}\,w(x)\quad\mbox{for every cube}\;Q\subseteq\mathbb R^n.$$
For the case $p=\infty$, $w\in A_\infty$, if for any given $\varepsilon>0$, we can find a positive number $\delta>0$ such that if $Q$ is a cube,
$E$ is a measurable subset of $Q$ with $|E|<\delta|Q|$, then $\int_E w(x)\,dx<\varepsilon\int_Q w(x)\,dx$.
\par
It is well known that $A_\infty=\underset{1<p<\infty}{\bigcup}A_p$, namely, a nonnegative, locally integrable function $w(x)$ satisfies the condition $A_\infty$ if and only if it satisfies the condition $A_p$ for some $1<p<\infty$. We also know that if $w\in A_p$ with $1<p<\infty$, then $w\in A_r$ for all $r>p$, and $w\in A_q$ for some $1<q<p$. Therefore, we will use the notation $q_w=\inf\{q>1:w\in A_q\}$ to denote the critical index of $w$. Obviously, if $w\in A_q$, $q>1$, then we have $1\le q_w<q$.
\par
Given a cube $Q$ and $\lambda>0$, $\lambda Q$ denotes the cube with the same center as $Q$ whose side length is $\lambda$ times that of $Q$. $Q=Q(x_0,r)$ denotes the cube centered at $x_0$ with side length $r$. For a weight function $w$ and a measurable set $E$, we set the weighted measure $w(E)=\int_E w(x)\,dx$, and we denote the characteristic function of $E$ by $\chi_{_E}$.
\par
We shall need the following lemmas. For the proofs of these results, please refer to [3, Chap \uppercase\expandafter{\romannumeral 4}] and [4, Chap 9].
\newtheorem*{lemmaA}{Lemma A}
\begin{lemmaA}
Let $w\in A_p$, $p\ge1$. Then, for any cube $Q$, there exists an absolute constant $C>0$ such that
$$w(2Q)\le C w(Q).$$
In general, for any $\lambda>1$, we have
$$w(\lambda Q)\le C\lambda^{np}w(Q),$$
where $C$ does not depend on $Q$ nor on $\lambda$.
\end{lemmaA}
\newtheorem*{lemmaB}{Lemma B}
\begin{lemmaB}
Let $w\in A_q$, $q>1$. Then, for all $r>0$, there exists a constant $C$ independent of $r$ such that
$$\int_{|x|\ge r}\frac{w(x)}{|x|^{nq}}\,dx\le C r^{-nq}w(Q(0,2r)).$$
\end{lemmaB}
\newtheorem*{lemmaC}{Lemma C}
\begin{lemmaC}
Let $w\in A_\infty$. For any $0<\varepsilon<1$, there exists a positive number $0<\delta<1$ such that if $E$ is a measurable subset of a cube $Q$ with $|E|/|Q|>\varepsilon$, then we have $w(E)/w(Q)>\delta.$
\end{lemmaC}
\newtheorem*{lemmaD}{Lemma D}
\begin{lemmaD}
Let $w\in A_p$, $p\ge1$. Then there exists an absolute constant $C>0$ such that $$C\left(\frac{|E|}{|Q|}\right)^p\le\frac{w(E)}{w(Q)},$$
for any measurable subset $E$ of a cube $Q$.
\end{lemmaD}
\par
Given a Muckenhoupt's weight function $w$ on $\mathbb R^n$, for $0<q<\infty$, we denote by $L^q_w(\mathbb R^n)$ the space of all functions satisfying
$$\|f\|_{L^q_w(\mathbb R^n)}=\left(\int_{\mathbb R^n}|f(x)|^qw(x)\,dx\right)^{1/q}<\infty.$$
When $q=\infty$, $L^\infty_w$ will be taken to mean $L^\infty$, and we set $\|f\|_{L^\infty_w}=\|f\|_{L^\infty}.$
As we all know, for any $0<p<\infty$, the weighted Hardy spaces $H^p_w(\mathbb R^n)$ can be defined in terms of maximal functions.
Let $\varphi$ be a function in $\mathscr S(\mathbb R^n)$ satisfying $\int_{\mathbb R^n}\varphi(x)\,dx=1$.
Set
$$\varphi_t(x)=t^{-n}\varphi(x/t),\quad t>0,\;x\in\mathbb R^n.$$
We will define the maximal function $M_\varphi f(x)$ by
$$M_\varphi f(x)=\sup_{t>0}|f*\varphi_t(x)|.$$
Then $H^p_w(\mathbb R^n)$ consists of those tempered distributions $f\in\mathscr S'(\mathbb R^n)$ for which
$M_\varphi f\in L^p_w(\mathbb R^n)$ with $\|f\|_{H^p_w}=\|M_\varphi f\|_{L^p_w}$. For every $1<p<\infty$, as in the unweighted case, we have $L^p_w(\mathbb R^n)=H^p_w(\mathbb R^n).$
\par
The real-variable theory of weighted Hardy spaces have been studied by many authors. In 1979, Garcia-Cuerva studied the atomic decomposition and the dual spaces of $H^p_w$ for $0<p\le1$. In 2002, Lee and Lin gave the molecular characterization of $H^p_w$ for $0<p\le1$, they also obtained the $H^p_w(\mathbb R)$, $\frac12<p\le1$ boundedness of the Hilbert transform and the $H^p_w(\mathbb R^n)$, $\frac n{n+1}<p\le1$ boundedness of the Riesz transforms. For the results mentioned above, we refer the readers to [2,6,9] for further details.
\par
In this article, we will use Garcia-Cuerva's atomic decomposition theory for weighted Hardy spaces in [2,9]. We characterize weighted Hardy spaces in terms of atoms in the following way.
\par
Let $0<p\le1\le q\le\infty$ and $p\ne q$ such that $w\in A_q$ with critical index $q_w$. Set [\,$\cdot$\,] the greatest integer function. For $s\in \mathbb Z_+$ satisfying $s\ge[n({q_w}/p-1)],$ a real-valued function $a(x)$ is called ($p,q,s$)-atom centered at $x_0$ with respect to $w$(or $w$-($p,q,s$)-atom centered at $x_0$) if the following conditions are satisfied:
\par (a) $a\in L^q_w(\mathbb R^n)$ and is supported in a cube $Q$ centered at $x_0$,
\par (b) $\|a\|_{L^q_w}\le w(Q)^{1/q-1/p}$,
\par (c) $\int_{\mathbb R^n}a(x)x^\alpha\,dx=0$ for every multi-index $\alpha$ with $|\alpha|\le s$.
\newtheorem*{theorem}{Theorem E}
\begin{theorem}
Let $0<p\le1\le q\le\infty$ and $p\ne q$ such that $w\in A_q$ with critical index $q_w$. For each $f\in H^p_w(\mathbb R^n)$, there exist a sequence \{$a_j$\} of $w$-$(p,q,[n(q_w/p-1)])$-atoms and a sequence \{$\lambda_j$\} of real numbers with $\sum_j|\lambda_j|^p\le C\|f\|^p_{H^p_w}$ such that $f=\sum_j\lambda_j a_j$ both in the sense of distributions and in the $H^p_w$ norm.
\end{theorem}
\noindent\textbf{\large{2. The intrinsic square functions and our main results}}
\par
The intrinsic square functions were first introduced by Wilson in [10] and [11],
the so-called intrinsic square functions are defined as follows.
For $0<\alpha\le1$, let ${\mathcal C}_\alpha$ be the family of functions $\varphi$ defined on $\mathbb R^n$ such that $\varphi$ has support containing in $\{x\in\mathbb R^n: |x|\le1\}$, $\int_{\mathbb R^n}\varphi(x)\,dx=0$ and for all $x, x'\in \mathbb R^n$,
$$|\varphi(x)-\varphi(x')|\le|x-x'|^\alpha.$$
For $(y,t)\in {\mathbb R}^{n+1}_+ =\mathbb R^n\times(0,\infty)$ and $f\in L^1_{{loc}}(\mathbb R^n)$, we set
$$A_\alpha(f)(y,t)=\sup_{\varphi\in{\mathcal C}_\alpha}|f*\varphi_t(y)|.$$
Then we define the intrinsic square function of $f$(of order $\alpha$) by the formula
$$S_\alpha(f)(x)=\left(\iint_{\Gamma(x)}\Big(A_\alpha(f)(y,t)\Big)^2\frac{dydt}{t^{n+1}}\right)^{1/2},$$
where $\Gamma(x)$ denotes the usual cone of aperture one:
$$\Gamma(x)=\{(y,t)\in{\mathbb R}^{n+1}_+:|x-y|<t\}.$$
We can also define varying-aperture versions of $S_\alpha(f)$ by the formula
$$S_{\alpha,\beta}(f)(x)=\left(\iint_{\Gamma_\beta(x)}\Big(A_\alpha(f)(y,t)\Big)^2\frac{dydt}{t^{n+1}}\right)^{1/2},$$
where $\Gamma_\beta(x)$ is the usual cone of aperture $\beta>0$:
$$\Gamma_\beta(x)=\{(y,t)\in{\mathbb R}^{n+1}_+:|x-y|<\beta t\}.$$
The intrinsic Littlewood-Paley $g$-function(could be viewed as ``zero-aperture" version of $S_\alpha(f)$) and the intrinsic $g^*_\lambda$-function(could be viewed as ``infinite aperture" version of $S_\alpha(f)$) will be defined respectively by
$$g_\alpha(f)(x)=\left(\int_0^\infty\Big(A_\alpha(f)(x,t)\Big)^2\frac{dt}{t}\right)^{1/2}$$
and
$$g^*_{\lambda,\alpha}(f)(x)=\left(\iint_{{\mathbb R}^{n+1}_+}\left(\frac t{t+|x-y|}\right)^{\lambda n}\Big(A_\alpha(f)(y,t)\Big)^2\frac{dydt}{t^{n+1}}\right)^{1/2}.$$
\par
Similarly, we can also introduce the so-called similar-looking square functions ${\tilde S}_{(\alpha,\varepsilon)}(f)(x)$, which are defined via convolutions with kernels that have unbounded supports, more precisely, for $0<\alpha\le 1$ and $\varepsilon>0$, let ${\mathcal C}_{(\alpha,\varepsilon)}$ be the family of functions $\varphi$ defined on $\mathbb R^n$ such that for all $x\in\mathbb R^n$,$$|\varphi(x)|\le(1+|x|)^{-n-\varepsilon},$$
and for all $x, x'\in\mathbb R^n$,
$$|\varphi(x)-\varphi(x')|\le|x-x'|^\alpha\left((1+|x|)^{-n-\varepsilon}+(1+|x'|)^{-n-\varepsilon}\right),$$
and also satisfy $\int_{\mathbb R^n}\varphi(x)\,dx=0$.
\\
Let $f$ be such that $|f(x)|(1+|x|)^{-n-\varepsilon}\in L^1(\mathbb R^n)$. For any $(y,t)\in {\mathbb R}^{n+1}_+$, set
$${\tilde A}_{(\alpha,\varepsilon)}(f)(y,t)=\sup_{\varphi\in {\mathcal C}_{(\alpha,\varepsilon)}}|f*\varphi_t(y)|.$$
We define
$${\tilde S}_{(\alpha,\varepsilon)}(f)(x)=\left(\iint_{\Gamma(x)}\Big({\tilde A}_{(\alpha,\varepsilon)}(f)(y,t)\Big)^2\frac{dydt}{t^{n+1}}\right)^{1/2},$$
$${\tilde g}_{(\alpha,\varepsilon)}(f)(x)=\left(\int_0^\infty\Big({\tilde A}_{(\alpha,\varepsilon)}(f)(x,t)\Big)^2\frac{dt}{t}\right)^{1/2},$$
and
$${\tilde g}^*_{\lambda,(\alpha,\varepsilon)}(f)(x)=\left(\iint_{{\mathbb R}^{n+1}_+}\left(\frac t{t+|x-y|}\right)^{\lambda n}\Big({\tilde A}_{(\alpha,\varepsilon)}(f)(y,t)\Big)^2\frac{dydt}{t^{n+1}}\right)^{1/2}.$$
\par
In [11], Wilson proved that the intrinsic square functions are bounded operators on the weighted Lebesgue spaces $L^p_w(\mathbb R^n)$ for $1<p<\infty$, namely, he showed the following result.
\newtheorem*{Theorem}{Theorem F}
\begin{Theorem}
Let $w\in A_p$, $1<p<\infty$ and $0<\alpha\le1$. Then there exists a positive constant $C>0$ such that$$\|S_\alpha(f)\|_{L^p_w}\le C \|f\|_{L^p_w}.$$
\end{Theorem}
Recently, Huang and Liu [5] studied the boundedness of intrinsic square functions on the weighted Hardy spaces $H^1_w(\mathbb R^n)$. Moreover, they obtained the intrinsic square function characterizations of $H^1_w(\mathbb R^n).$
\par
As a continuation of their work, the purpose of this paper is to investigate the boundedness of intrinsic square functions on the weighted Hardy spaces $H^p_w(\mathbb R^n)$ for $0<p<1$. Furthermore, we will characterize the weighted Hardy spaces $H^p_w(\mathbb R^n)$ for $0<p<1$ by the intrinsic square functions including the Lusin area function, Littlewood-Paley $g$-function and $g^*_\lambda$-function.
\par
In order to state our theorems, we need to introduce the Lipschitz space $Lip(\alpha,1,0)$ for $0<\alpha\le1$.
Set $b_Q=\frac{1}{|Q|}\int_Q b(x)\,dx$.
$$Lip(\alpha,1,0)=\{b\in L_{loc}(\mathbb R^n):\|b\|_{Lip(\alpha,1,0)}<\infty\},$$
where
$$\|b\|_{Lip(\alpha,1,0)}=\sup_Q\frac{1}{|Q|^{1+\alpha/n}}\int_Q|b(y)-b_Q|\,dy$$
and the supremum is taken over all cubes $Q$ in $\mathbb R^n$.
\par
Our main results are stated as follows.
\newtheorem{theorem1}{Theorem}
\begin{theorem1}
Let $0<\alpha\le1$, $\frac{n}{n+\alpha}<p<1$, $w\in A_{p(1+\frac\alpha n)}$ and $\varepsilon>\alpha$. Suppose that $f\in \big(Lip(\alpha,1,0)\big)^*$, then a tempered distribution $f\in H^p_w(\mathbb R^n)$ if and only if $g_\alpha(f)\in L^p_w(\mathbb R^n)$ or ${\tilde g}_{(\alpha,\varepsilon)}(f)\in L^p_w(\mathbb R^n)$ and $f$ vanishes weakly at infinity.
\end{theorem1}
\begin{theorem1}
Let $0<\alpha\le1$, $\frac{n}{n+\alpha}<p<1$, $w\in A_{p(1+\frac\alpha n)}$ and $\varepsilon>\alpha$. Suppose that $f\in \big(Lip(\alpha,1,0)\big)^*$, then a tempered distribution $f\in H^p_w(\mathbb R^n)$ if and only if $S_\alpha(f)\in L^p_w(\mathbb R^n)$ or ${\tilde S}_{(\alpha,\varepsilon)}(f)\in L^p_w(\mathbb R^n)$ and $f$ vanishes weakly at infinity.
\end{theorem1}
\begin{theorem1}
Let $0<\alpha\le1$, $\frac{n}{n+\alpha}<p<1$, $w\in A_{p(1+\frac\alpha n)}$, $\varepsilon>\alpha$ and $\lambda>\frac{3n+2\alpha}{n}$. Suppose that $f\in \big(Lip(\alpha,1,0)\big)^*$, then a tempered distribution $f\in H^p_w(\mathbb R^n)$ if and only if $g^*_{\lambda,\alpha}(f)\in L^p_w(\mathbb R^n)$ or ${\tilde g}^*_{\lambda, (\alpha,\varepsilon)}(f)\in L^p_w(\mathbb R^n)$ and $f$ vanishes weakly at infinity.
\end{theorem1}
\noindent\textbf{Remark 1.} Clearly, if for every $t>0$, $\varphi_t\in \mathcal C _\alpha$, then we have $\varphi_t\in Lip(\alpha,1,0)$. Thus the intrinsic square functions are well defined for tempered distributions in $\big(Lip(\alpha,1,0)\big)^*$.\\
\textbf{Remark 2.} We say that a tempered distribution $f$ vanishes weakly at infinity, if for any $\varphi\in\mathscr S$, we have $f*\varphi_t(x)\to0$ as $t\to\infty$ in the sense of distributions.
\par
Throughout this article, we will use $C$ to denote a positive constant, which is independent of the main parameters and not necessarily the same at each occurrence. By $A\sim B$, we mean that there exists a constant $C>1$ such that $\frac1C\le\frac AB\le C$. \\
\textbf{\large{3. The necessity of our theorems}}
\par
We shall first prove the following lemma.
\newtheorem*{lemma3}{Lemma 3.1}
\begin{lemma3}
Let $0<p<1$ and $w\in A_\infty$. Then for every $f\in H^p_w(\mathbb R^n)$, we have that $f$ vanishes weakly at infinity.
\end{lemma3}
\begin{proof}
For any given $\varphi\in\mathscr S(\mathbb R^n)$, $\int_{\mathbb R^n}\varphi(x)\,dx=1$, we denote the nontangential maximal function of $f$ by
$$M^*_{\varphi}(f)(x)=\sup_{|y-x|<t}|f*\varphi_t(y)|.$$
Then we have $|f*\varphi_t(x)|\le M^*_{\varphi}(f)(y)$ whenever $|x-y|<t$. As a consequence, we have the following inequality
$$\int_{|x-y|<t}|f*\varphi_t(x)|^pw(y)\,dy\le\int_{|x-y|<t}\left(M^*_{\varphi}(f)(y)\right)^pw(y)\,dy.$$
Hence
$$|f*\varphi_t(x)|^p\le\frac{1}{w(Q(x,\sqrt2 t))}\|M^*_{\varphi}(f)\|_{L^p_w}^p\le C\frac{1}{w(Q(x,\sqrt2 t))}\|M_{\varphi}(f)\|_{L^p_w}^p.$$
It is well known that for given $w\in A_\infty$, $w$ satisfies the doubling condition(Lemma A). Furthermore, we can easily prove that $w$ also satisfies the reverse doubling condition; that is, for any cube $Q$, there exists a constant $C_1>1$ such that $w(2Q)\ge C_1 w(Q)$. From this property, we can deduce $w(2^k Q)\ge C_1^kw(Q)$ by induction. Set $Q=Q(x,\sqrt 2)$. So we can get
$$\lim_{k\to\infty}\frac{1}{w(2^k Q)}=0,$$
which implies
$$\lim_{t\to\infty}\frac{1}{w(Q(x,\sqrt2 t))}=0.$$
This completes the proof of the lemma.
\end{proof}
\par
From the definitions of intrinsic square functions, we know that when $\varphi\in{\mathcal C}_\alpha$, $0<\alpha\le1$, then there exists a positive constant $c$ depending only on $\alpha,\varepsilon,$ and $n$, such that $c\varphi\in{\mathcal C}_{(\alpha,\varepsilon)}$. Thus we can get the pointwise inequality $S_\alpha(f)(x)\le C{\tilde S}_{(\alpha,\varepsilon)}(f)(x)$. Furthermore, in [10], the author proved that this inequality have a partial converse; that is, for every $\alpha'$ satisfying $0<\alpha'\le\alpha$ and $\alpha'<\varepsilon$, for all $f$ such that $|f(x)|(1+|x|)^{-n-\varepsilon}\in L^1(\mathbb R^n)$, we have ${\tilde S}_{(\alpha,\varepsilon)}(f)(x)\le C S_\alpha(f)(x)$. So if we choose $\alpha'=\alpha$ and $\varepsilon>\alpha$, we obtain $S_\alpha(f)(x)\sim{\tilde S}_{(\alpha,\varepsilon)}(f)(x).$ In [10], the author also showed that the functions $S_\alpha(f)(x)$ and $g_\alpha(f)(x)$ are pointwise comparable. Meanwhile, he pointed out that by similar arguments we can show the pointwise comparability of ${\tilde S}_{(\alpha,\varepsilon)}(f)(x)$ and ${\tilde g}_{(\alpha,\varepsilon)}(f)(x)$. Therefore, in order to prove the necessity of Theorems 1, 2, we need only to prove the following proposition.
\newtheorem*{prop}{Proposition 3.2}
\begin{prop}
Let $0<\alpha\le1$, $\frac{n}{n+\alpha}<p<1$ and $w\in A_{p(1+\frac\alpha n)}$. Then for every $f\in H^p_w(\mathbb R^n)$, we have
$$\|g_\alpha(f)\|_{L^p_w}\le C\|f\|_{H^p_w}.$$
\end{prop}
\begin{proof}
Set $q=p(1+\frac\alpha n)$. Then for $w\in A_q$, we have $[n({q_w}/p-1)]=0.$ By Theorem E, it suffices to show that for any $w$-$(p,q,0)$-atom $a$, there exists a constant $C>0$ independent of $a$ such that $\|g_\alpha(a)\|_{L^p_w}\le C.$
\par
Let $a$ be a $w$-$(p,q,0)$-atom with supp $a\subset Q=Q(x_0,r)$, and let $Q^*=2\sqrt nQ.$ By using H\"{o}lder inequality, Lemma A and Theorem F, we have
\begin{align}
\int_{Q^*}|g_\alpha(a)(x)|^pw(x)\,dx&\le\left(\int_{Q^*}|g_\alpha(a)(x)|^qw(x)\,dx\right)^{p/q}\left(\int_{Q^*}w(x)\,dx\right)^{1-p/q}\notag\\
                                 &\le \|g_\alpha(a)\|^p_{L^q_w}w(Q^*)^{1-p/q}\notag\\
                                 &\le C\|S_\alpha(a)\|^p_{L^q_w}w(Q)^{1-p/q}\notag\\
                                 &\le C\|a\|^p_{L^q_w}w(Q)^{1-p/q}\\
                                 &\le C\notag.
\end{align}
Below we give the estimate of the integral $I=\int_{(Q^*)^c}|g_\alpha(a)(x)|^pw(x)\,dx.$\par
For any $\varphi\in{\mathcal C}_\alpha$, by the vanishing moment condition of atom $a$, we have
\begin{equation}
\begin{split}
\big|a*\varphi_t(x)\big|&=\left|\int_Q\big(\varphi_t(x-y)-\varphi_t(x-x_0)\big)a(y)\,dy\right|\\
&\le\int_Q\frac{|y-x_0|^\alpha}{t^{n+\alpha}}|a(y)|\,dy\\
&\le C\cdot\frac{r^\alpha}{t^{n+\alpha}}\int_Q|a(y)|\,dy.
\end{split}
\end{equation}
Denote the conjugate exponent of $q>1$ by $q'=q/{(q-1)}$. H\"{o}lder's inequality and the $A_q$ condition yield
\begin{equation}
\begin{split}
\int_Q|a(y)|\,dy&\le\left(\int_Q|a(y)|^qw(y)\,dy\right)^{1/q}\left(\int_Q w(y)^{-1/(q-1)}\,dy\right)^{1/q'}\\
&\le C\|a\|_{L^q_w}\left(\frac{|Q|^q}{w(Q)}\right)^{1/q}\\
&\le C\frac{|Q|}{w(Q)^{1/p}}.
\end{split}
\end{equation}
We note that supp $\varphi\subseteq\{x\in\mathbb R^n:|x|\le1\}$, then for any $y\in Q$, $x\in(Q^*)^c$, we have $t\ge|x-y|\ge|x-x_0|-|y-x_0|\ge\frac{|x-x_0|}{2}.$
\\Substituting the above inequality (3) into (2), we thus obtain
\begin{equation}
\begin{split}
|g_\alpha(a)(x)|^2&=\int_0^\infty\Big(\sup_{\varphi\in{\mathcal C}_\alpha}\big|a*\varphi_t(x)\big|\Big)^2\frac{dt}{t}\\
&\le C\left(\frac{|Q|}{w(Q)^{1/p}}\right)^2r^{2\alpha}\int_{\frac{|x-x_0|}{2}}^\infty\frac{dt}{t^{2(n+\alpha)+1}}\\
&\le C\left(\frac{|Q|}{w(Q)^{1/p}}\right)^2r^{2\alpha}\frac{1}{|x-x_0|^{2n+2\alpha}}
\end{split}
\end{equation}
It follows from (4), Lemma A and Lemma B that
\begin{equation}
\begin{split}
I&=\int_{(Q^*)^c}|g_\alpha(a)(x)|^pw(x)\,dx\\
  &\le C\left(\frac{r^{n+\alpha}}{w(Q)^{1/p}}\right)^p\int_{|x-x_0|\ge\sqrt n r}\frac{w(x)}{|x-x_0|^{nq}}\,dx\\
  &=C\left(\frac{r^{n+\alpha}}{w(Q)^{1/p}}\right)^p\int_{|y|\ge\sqrt n r}\frac{w_1(y)}{|y|^{nq}}\,dy\\
  &\le C\left(\frac{r^{n+\alpha}}{w(Q)^{1/p}}\right)^p r^{-nq}w_1(Q_1)\\
  &=C\left(\frac{r^{n+\alpha}}{w(Q)^{1/p}}\right)^p r^{-nq}w(Q)\\
  &\le C,
\end{split}
\end{equation}
where $w_1(x)=w(x+x_0)$ is the translation of $w(x)$, $Q_1$ is a cube which is the translation of $Q$. It is obvious that $w_1\in A_q$ for $w\in A_q$, $q>1$, and $q_{w_1}=q_w$.
Therefore, Proposition 3.2 is proved by combining (1) and (5).
\end{proof}
\newtheorem*{Prop}{Proposition 3.3}
\begin{Prop}
Let $0<\alpha\le1$, $\frac{n}{n+\alpha}<p<1$, $w\in A_{p(1+\frac\alpha n)}$ and $\lambda>\frac{3n+2\alpha}{n}$. Then for every $f\in H^p_w(\mathbb R^n)$, we have$$\|g^*_{\lambda,\alpha}(f)\|_{L^p_w}\le C\|f\|_{H^p_w}.$$
\end{Prop}
\begin{proof}
Let $q=p(1+\frac{\alpha}{n})$. As in the proof of Proposition 3.2, we only need to show that for any $w$-$(p,q,0)$-atom $a$, there exists a constant $C>0$ independent of $a$ such that $\|g^*_{\lambda,\alpha}(a)\|_{L^p_w}\le C.$
\par
Let $a$ be a $w$-$(p,q,0)$-atom with supp $a\subset Q=Q(x_0,r)$, and let $Q^*_k=2\sqrt n(2^kQ)$. From the definition, we readily see that
\begin{equation*}
\begin{split}
g^*_{\lambda,\alpha}(a)(x)^2&=\iint_{\mathbb R^{n+1}_+}\left(\frac{t}{t+|x-y|}\right)^{\lambda n}\Big(A_\alpha(a)(y,t)\Big)^2\frac{dydt}{t^{n+1}}\\
&=\int_0^\infty\int_{|x-y|<t}\left(\frac{t}{t+|x-y|}\right)^{\lambda n}\Big(A_\alpha(a)(y,t)\Big)^2\frac{dydt}{t^{n+1}}\\
&+\sum_{k=1}^\infty\int_0^\infty\int_{2^{k-1}t\le|x-y|<2^kt}\left(\frac{t}{t+|x-y|}\right)^{\lambda n}\Big(A_\alpha(a)(y,t)\Big)^2\frac{dydt}{t^{n+1}}\\
&\le C\bigg[S_\alpha(a)(x)^2+\sum_{k=1}^\infty 2^{-k\lambda n}S_{\alpha,2^k}(a)(x)^2\bigg].
\end{split}
\end{equation*}
Since $0<p<1$, we thus get
$$\|g^*_{\lambda,\alpha}(a)\|_{L^p_w}^p\le C\bigg[\|S_\alpha(a)\|^p_{L^p_w}+\sum_{k=1}^\infty 2^{-\frac{k\lambda np}{2}}\|S_{\alpha,2^k}(a)\|^p_{L^p_w}\bigg].$$
By Proposition 3.2, we can obtain $\|S_\alpha(a)\|_{L^p_w}\le C$.
It remains to estimate $\|S_{\alpha,2^k}(a)\|_{L^p_w}$ for $k=1,2,\ldots.$
\\
First we claim that the following inequality holds.
\begin{equation}\|S_{\alpha,2^k}(a)\|_{L^2_w}\le C\cdot2^{\frac{knq}{2}}\|S_\alpha(a)\|_{L^2_w} \quad k=1,2,\ldots.\end{equation}
In fact, by the Fubini theorem and Lemma A, we can get
\begin{equation*}
\begin{split}
\|S_{\alpha,2^k}(a)\|_{L^2_w}^2&=\int_{\mathbb R^n}\bigg(\int_{{\mathbb R}^{n+1}_+}\Big(A_\alpha(a)(y,t)\Big)^2\chi_{|x-y|<2^k t}\frac{dydt}{t^{n+1}}\bigg)w(x)\,dx\\
&=\int_{{\mathbb R}^{n+1}_+}\Big(\int_{|x-y|<2^k t}w(x)\,dx\Big)\Big(A_\alpha(a)(y,t)\Big)^2\frac{dydt}{t^{n+1}}\\
&\le C\cdot2^{knq}\int_{{\mathbb R}^{n+1}_+}\Big(\int_{|x-y|<t}w(x)\,dx\Big)\Big(A_\alpha(a)(y,t)\Big)^2\frac{dydt}{t^{n+1}}\\
&=C\cdot 2^{knq}\|S_\alpha(a)\|_{L^2_w}^2.
\end{split}
\end{equation*}
Using H\"{o}lder's inequality, Lemma A, Theorem F and (6), we obtain
\begin{equation}
\begin{split}
\Big(\int_{Q^*_k}|S_{\alpha,2^k}(a)(x)|^pw(x)\,dx\Big)^{1/p}&\le\|S_{\alpha,2^k}(a)\|_{L^2_w}w(Q^*_k)^{\frac 1p-\frac12}\\
&\le C\cdot2^{\frac{knq}{2}}\|S_\alpha(a)\|_{L^2_w}\big(2^{knq}w(Q)\big)^{\frac 1p-\frac12}\\
&\le C\cdot2^{\frac{knq}{p}}\|a\|_{L^2_w}\left(w(Q)\right)^{\frac 1p-\frac12}\\
&\le C\cdot2^{\frac{knq}{p}},
\end{split}
\end{equation}
where we have used the fact that $w\in A_q$, $1<q<1+\frac{\alpha}{n}\le2$, then $w\in A_2$.
\\
Below we give the estimate of the integral $J=\int_{(Q^*_k)^c}|S_{\alpha,2^k}(a)(x)|^pw(x)\,dx$.
\par
Note that supp $\varphi\subset\{x\in\mathbb R^n:|x|\le1\}$, by a simple calculation, we know that for any $(y,t)\in\Gamma_{2^k}(x)$, $x\in(Q^*_k)^c$, then $t\ge\frac{|x-x_0|}{2^{k+1}}$. It follows from (2) and (3) that
\begin{equation}
\begin{split}
\left|S_{\alpha,2^k}(a)(x)\right|^2&\le C\left(\frac{|Q|}{w(Q)^{1/p}}\right)^2r^{2\alpha}\iint_{\Gamma_{2^k}(x)}\frac{dydt}{t^{2(n+\alpha)}\cdot t^{n+1}}\\
&\le C\left(\frac{|Q|}{w(Q)^{1/p}}\right)^2r^{2\alpha}2^{kn}\int_{\frac{|x-x_0|}{2^{k+1}}}^\infty\frac{dt}{t^{2(n+\alpha)+1}}\\
&\le C\cdot2^{3kn+2k\alpha}\left(\frac{r^{n+\alpha}}{w(Q)^{1/p}}\right)^2\frac{1}{|x-x_0|^{2(n+\alpha)}}.
\end{split}
\end{equation}
Using Lemma A, Lemma B and (8), we have
\begin{equation}
\begin{split}
J&=\int_{(Q^*_k)^c}|S_{\alpha,2^k}(a)(x)|^pw(x)\,dx\\
&\le C\cdot2^{\frac{kp(3n+2\alpha)}{2}}\frac{r^{p(n+\alpha)}}{w(Q)}\int_{|x-x_0|\ge{\sqrt n}2^k r}\frac{w(x)}{|x-x_0|^{nq}}\,dx\\
&\le C\cdot2^{\frac{kp(3n+2\alpha)}{2}}\frac{r^{p(n+\alpha)}}{w(Q)}(2^k r)^{-nq}(2^k)^{nq}w_1(Q_1)\\
&\le C\cdot2^{\frac{kp(3n+2\alpha)}{2}},
\end{split}
\end{equation}
where the notations $w_1$ and $Q_1$ are the same as Proposition 3.2, we have $w_1(Q_1)=w(Q)$. Hence, by the above estimates (7) and (9), we obtain
$$\|S_{\alpha,2^k}(a)\|^p_{L^p_w}\le C\cdot\Big(2^{kp(n+\alpha)}+2^{\frac{kp(3n+2\alpha)}{2}}\Big)\le C\cdot2^{\frac{kp(3n+2\alpha)}{2}}.$$
Therefore
$$\|g^*_{\lambda,\alpha}(a)\|_{L^p_w}^p\le C\sum_{k=1}^\infty 2^{-\frac{k\lambda np}{2}}\cdot2^{\frac{kp(3n+2\alpha)}{2}}\le C,$$
where the last inequality holds since $\lambda>3+2\alpha/n.$
The proof of Proposition 3.3 is complete.
\end{proof}Using the same arguments as above, we can also show the $H^p_w$-$L^p_w$ boundedness of ${\tilde g}^*_{\lambda, (\alpha,\varepsilon)}$; that is,
\begin{equation}
\|{\tilde g}^*_{\lambda, (\alpha,\varepsilon)}(f)\|_{L^p_w}\le C\|f\|_{H^p_w}.
\end{equation}
Therefore, by Lemma 3.1, Proposition 3.2, Proposition 3.3 and (10), we have proved the necessity of Theorems 1, 2 and 3.\\
\textbf{\large{4. The sufficiency of our theorems}}\par
We shall need the following Calder\'{o}n reproducing formula given in [1].
\newtheorem*{lemma4}{Lemma 4.1}
\begin{lemma4}
Let $\psi\in\mathscr S(\mathbb R^n)$, supp ${\psi\subset\{x\in\mathbb R^n:|x|\le1\}}$, $\int_{\mathbb R^n}\psi(x)\,dx=0$ and
\begin{equation*}
\int_0^\infty\big|\hat\psi(\xi t)\big|^2\frac{dt}{t}=1 \quad \mbox{whenever}\; \;\xi\ne0.
\end{equation*}
Then for any $f\in\mathscr S'(\mathbb R^n)$,
$f$ vanishes weakly at infinity, we have
\begin{equation}
f(x)=\int_0^\infty\int_{\mathbb R^n}f*\psi_t(y)\psi_t(x-y)\frac{dydt}{t},
\end{equation}
where the equality holds in the sense of distribution.
\end{lemma4}
Suppose that $\psi$ satisfies the conditions of Lemma 4.1. For every $f\in\mathscr S'(\mathbb R^n)$, we define the area integral of $f$ by
$$S_\psi(f)(x)=\left(\int_{|x-y|<t}\big|f*\psi_t(y)\big|^2\frac{dydt}{t^{n+1}}\right)^{1/2}.$$
\par We now prove the following result.
\newtheorem*{prop4}{Proposition 4.2}
\begin{prop4}
Let $0<\alpha\le1$, $\frac{n}{n+\alpha}<p<1$ and $w\in A_{p(1+\frac{\alpha}{n})}$. Then for any $f\in\mathscr S'(\mathbb R^n)$, $f$ vanishes weakly at infinity, we have
$$\|f\|_{H^p_w}\le C\|S_\psi(f)\|_{L^p_w}.$$
\end{prop4}
\begin{proof}
We follow the same constructions as in [7]. For any $k\in\mathbb Z$, set
$$\Omega_k=\{x\in\mathbb R^n:S_\psi(f)(x)>2^k\}.$$
Let $\mathbb D$ denote the set formed by all dyadic cubes in $\mathbb R^n$ and let
$$\mathbb D_k=\Big\{Q\in{\mathbb D}:|Q\cap\Omega_k|>\frac{|Q|}{2},|Q\cap\Omega_{k+1}|\le\frac{|Q|}{2}\Big\}.$$
Obviously, for any $Q\in\mathbb D$, there exists a unique $k\in\mathbb Z$ such that $Q\in{\mathbb D}_k.$
We also denote the maximal dyadic cubes in ${\mathbb D}_k$ by $Q_k^l$. Set
$$\widetilde Q=\{(y,t)\in{\mathbb R}^{n+1}_+:y\in Q,\, l(Q)<t\le 2l(Q)\},$$
where $l(Q)$ denotes the side length of $Q$.\\
If we set $\widetilde{Q_k^l}=\underset{Q_k^l\supseteq Q\in{\mathbb D}_k}{\bigcup}\widetilde Q$,
then we have ${\mathbb R}^{n+1}_+=\underset{k}{\bigcup}\,\underset{l}{\bigcup}\,\widetilde{Q^l_k}$. Hence, by (11), we obtain
$$f(x)=\sum_k\sum_l\int_{\widetilde{Q^l_k}}f*\psi_t(y)\psi_t(x-y)\frac{dydt}{t}=\sum_k\sum_l\lambda_{kl}a^l_k(x),$$
where
$$a_k^l(x)=\lambda_{kl}^{-1}\int_{\widetilde{Q^l_k}}f*\psi_t(y)\psi_t(x-y)\frac{dydt}{t}$$
and
$$\lambda_{kl}=w(Q^l_k)^{1/p-1/2}\bigg(\int_{\widetilde{Q^l_k}}\big|f*\psi_t(y)\big|^2\frac{w(Q^l_k)}{|Q^l_k|}\frac{dydt}{t}\bigg)^{1/2}.$$
By the properties of $\psi$, we can easily get supp $a^l_k\subseteq5Q^l_k$, $\int_{\mathbb R^n}a^l_k(x)\,dx=0.$\\
Let $q=p(1+\frac{\alpha}{n})$, $w\in A_q$. Since
$$\|a^l_k\|_{_{L^q_w}}=\sup_{\|b\|_{L_w^{q'}}\le1}\Big|\int_{\mathbb R^n}a^l_k(x)b(x)w(x)\,dx\Big|.$$
Then H\"{o}lder's inequality and the definition of $\lambda_{kl}$ imply
\begin{equation*}
\begin{split}
&\Big|\int_{\mathbb R^n}a^l_k(x)b(x)w(x)\,dx\Big|\\
\le&\lambda_{kl}^{-1}\int_{\widetilde{Q^l_k}}\big|f*\psi_t(y)\big|\big|g*\psi_t(y)\big|\frac{dydt}{t}\\
\le&\lambda_{kl}^{-1}\bigg(\int_{\widetilde{Q^l_k}}\big|f*\psi_t(y)\big|^2\frac{dydt}{t}\bigg)^{1/2}\bigg(\int_{\widetilde{Q^l_k}}\big|g*\psi_t(y)\big|^2\frac{dydt}{t}\bigg)^{1/2}\\
\le&\frac{|Q^l_k|^{1/2}}{w(Q^l_k)^{1/p}}\bigg(\int_{\widetilde{Q^l_k}}\big|g*\psi_t(y)\big|^2\frac{dydt}{t}\bigg)^{1/2},
\end{split}
\end{equation*}
where $g(x)=\chi_{5Q^l_k}(x)b(x)w(x)$. A simple calculation shows that
$$|g*\psi_t(y)|\le C\cdot{t^{-n}}\|b\|_{L_w^{q'}}w(Q^l_k)^{1/q}.$$
Hence we have
\begin{equation*}
\begin{split}
\|a^l_k\|_{L^q_w}&\le C\cdot\frac{|Q^l_k|^{1/2}}{w(Q^l_k)^{1/p}}w(Q^l_k)^{1/q}\bigg(\int_{\widetilde{Q^l_k}}\frac{dydt}{t^{2n+1}}\bigg)^{1/2}\\
&\le C\cdot w(Q^l_k)^{1/q-1/p},
\end{split}
\end{equation*}
where in the last inequality we have used the fact that for any $(y,t)\in\widetilde{Q^l_k}$, we have $t^n\sim|Q^l_k|$. Therefore these functions $a_k^l$ defined above are all $w$-$(p,q,0)$-atoms.\par
Set
$\Omega_k^*=\Big\{x\in\mathbb R^n:M_w(\chi_{_{\Omega_k}})(x)>\frac{C_0}{2}\Big\},$ where $C_0$ is a constant to be determined later and
$M_w(f)(x)=\sup_{x\in Q}\frac{1}{w(Q)}\int_Q|f(y)|w(y)\,dy$. Using the weighted weak type estimate of weighted maximal operator $M_w$, we have
$w(\Omega_k^*)\le Cw(\Omega_k)$. Consequently
$$\int_{\Omega^*_k \backslash {\Omega_{k+1}}}S_\psi(f)(x)^2w(x)\,dx\le(2^{k+1})^2w(\Omega_k^*)\le C\cdot2^{2k}w(\Omega_k).$$
We set $E=\{x\in\Omega^*_k\backslash\Omega_{k+1}:|x-y|<t\}$, then we have
\begin{equation*}
\begin{split}
\int_{\Omega^*_k\backslash\Omega_{k+1}}S_\psi(f)(x)^2w(x)\,dx&=\int_{{\mathbb R}^{n+1}_+}\Big\{\int_{\mathbb R^n}\chi_E(x)w(x)\,dx\Big\}\big|f*\psi_t(y)\big|^2\frac{dydt}{t^{n+1}}\\
&\ge\sum_{Q\in{\mathbb D}_k}\left(\int_{\widetilde Q}\big|f*\psi_t(y)\big|^2\frac{dydt}{t^{n+1}}\right)w(E).
\end{split}
\end{equation*}
For any $x\in Q$, $Q\in{\mathbb D}_k$, we have $|Q\cap\Omega_k|>\frac12|Q|$, which is equivalent to
\begin{equation}
\frac{1}{|Q|}\int_Q\chi_{_{\Omega_k}}(y)\,dy>\frac12.
\end{equation}
H\"{o}lder's inequality and the $A_q$ condition give
\begin{align}
\frac{1}{|Q|}\int_Q\chi_{_{\Omega_k}}(y)\,dy&\le\frac{1}{|Q|}\left(\int_Q|\chi_{_{\Omega_k}}(y)|^qw(y)\,dy\right)^{1/q}\left(\int_Q w^{-1/(q-1)}\,dy\right)^{(q-1)/q}\notag\\
&\le[w]_{A_q}^{\frac1q}\left(\frac{1}{w(Q)}\int_Q\chi_{_{\Omega_k}}(y)w(y)\,dy\right)^{1/q}.
\end{align}
It follows immediately from (12) and (13) that $M_w(\chi_{_{\Omega_k}})(x)>\left(\frac12\right)^q[w]_{A_q}^{-1}$.
So if we choose $C_0=\left(\frac12\right)^{q-1}[w]_{A_q}^{-1}$, we have $x\in\Omega^*_k,$ which implies $Q\subseteq\Omega^*_k$. Hence $w(Q\cap\Omega^*_k)=w(Q).$
Since $|Q\cap\Omega_{k+1}|\le\frac12|Q|$, $w\in A_\infty$, then there exists a constant $0<C'<1$ such that $w(Q\cap\Omega_{k+1})\le C'w(Q)$. Consequently
\begin{equation}
\begin{split}
w(E)&\ge w(Q\cap(\Omega^*_k\backslash\Omega_{k+1}))\\
&\ge w(Q)-w(Q\cap\Omega_{k+1})\\
&\ge(1-C')w(Q).
\end{split}
\end{equation}
Suppose that $Q^l_k$ is the maximal dyadic cubes containing $Q$ which belong to ${\mathbb D}_k$. Then by Lemma D and the inequality (14), we can get
\begin{equation}
\begin{split}
2^{2k}w(\Omega_k)&\ge C\sum_{Q\in{\mathbb D}_k}\int_{\widetilde Q}\big|f*\psi_t(y)\big|^2w(Q)\frac{dydt}{t^{n+1}}\\
&\ge C\sum_{Q\in{\mathbb D}_k}\int_{\widetilde Q}\big|f*\psi_t(y)\big|^2w(Q^l_k)\left(\frac{|Q|}{|Q^l_k|}\right)^q\frac{dydt}{t^{n+1}}\\
&\ge C\sum_l\int_{\widetilde{Q^l_k}}\big|f*\psi_t(y)\big|^2\frac{w(Q^l_k)}{|Q^l_k|}\cdot\frac{1}{|Q^l_k|^{\frac\alpha n}}\frac{dydt}{t^{1-\alpha}}\\
&\ge C\sum_l\int_{\widetilde{Q^l_k}}\big|f*\psi_t(y)\big|^2\frac{w(Q^l_k)}{|Q^l_k|}\frac{dydt}{t},
\end{split}
\end{equation}
where the last inequality holds since $t\sim l(Q^l_k)$. For any $l\in{\mathbb Z}_+$, since $|Q^l_k\cap\Omega_k|>\frac12|Q^l_k|$, $w\in A_\infty$, then by Lemma C, we have that there exists a constant $0<C''<1$ such that $w(Q_k^l\cap\Omega_k)>C''w(Q^l_k).$ Note that the maximal dyadic cubes $Q^l_k$ are pairwise disjoint, we thus obtain
\begin{equation}
\begin{split}
w(\Omega_k)&\ge w\Big(\big(\underset{l}{\cup} Q^l_k\big)\cap\Omega_k\Big)\\
&=\sum_lw(Q_k^l\cap\Omega_k)\\
&>C''\sum_l w(Q^l_k).
\end{split}
\end{equation}
Then it follows from H\"{o}lder's inequality, (15) and (16) that
\begin{equation*}
\begin{split}
\sum_k\sum_l|\lambda_{kl}|^p&=\sum_k\sum_l\Big(w(Q^l_k)\Big)^{1-p/2}\bigg(\int_{\widetilde{Q^l_k}}\big|f*\psi_t(y)\big|^2\frac{w(Q^l_k)}{|Q^l_k|}\frac{dydt}{t}\bigg)^{p/2}\\
&\le\sum_k\Big(\sum_lw(Q^l_k)\Big)^{1-p/2}\bigg(\sum_l\int_{\widetilde{Q^l_k}}\big|f*\psi_t(y)\big|^2\frac{w(Q^l_k)}{|Q^l_k|}\frac{dydt}{t}\bigg)^{p/2}
\end{split}
\end{equation*}
\begin{equation*}
\begin{split}
&\le C\sum_k\Big(w(\Omega_k)\Big)^{1-p/2}\Big(2^{2k}w(\Omega_k)\Big)^{p/2}\\
&\le C\|S_\psi(f)\|^p_{L^p_w}.
\end{split}
\end{equation*}
Therefore, by using the atomic decomposition of weighted Hardy spaces, we get the desired result.
\end{proof}
Finally, we choose a function $\psi$ satisfying the conditions of Lemma 4.1. Obviously, we have
$\psi\in{\mathcal C}_\alpha$ for $0<\alpha\le1$, which implies
\begin{equation}
S_\psi(f)(x)\le S_\alpha(f)(x)\le C{\tilde S}_{(\alpha,\varepsilon)}(f)(x)\le C{\tilde g}^*_{\lambda,(\alpha,\varepsilon)}(f)(x).
\end{equation}
Combining the above inequality (17) and Proposition 4.2, we have proved the sufficiency of Theorems 1, 2 and 3.

\end{document}